%% file: sample.tex
\documentclass[letterpaper, 10 pt, conference]{ieeeconf}  %

\usepackage[utf8]{inputenc}

\IEEEoverridecommandlockouts                              %
\usepackage{graphics} %
\usepackage{epsfig} %
\usepackage{amsmath} %
\usepackage{amssymb} %
\usepackage{cite}
\usepackage{graphicx}
\usepackage{subcaption}

\title{\LARGE \bf
Reduced Order Modeling for Nonlinear PDE-constrained Optimization using Neural Networks
}

\author{Nikolaj Takata Mücke, Lasse Hjuler Christiansen, Allan Peter Engsig-Karup, John Bagterp J{\o}rgensen%
\thanks{Nikolaj Takata Mücke, Lasse Hjuler Christiansen, Allan Peter Karup-Engsig and John Bagterp Jørgensen are with the Department of Applied Mathematics and Computer Science, Technical University of Denmark,
	DK-2800 Kgs. Lyngby, Denmark \{ntmy,lhch,apek,jbjo\}@dtu.dk}
}

\begin{document}
\maketitle
\thispagestyle{empty}
\pagestyle{empty}

\begin{abstract}
\input{abstract}
\end{abstract}

\input{intro}

\input{theory}

\input{ex}

\input{concl}

\bibliographystyle{IEEEtran}
\bibliography{IEEEabrv,References}

\end{document}

%% file: abstract.tex
Nonlinear model predictive control (NMPC) often requires real-time solution to optimization problems. However, in cases where the mathematical model is of high dimension in the solution space, e.g. for solution of partial differential equations (PDEs), black-box optimizers are rarely sufficient to get the required online computational speed. In such cases one must resort to customized solvers. This paper present a new solver for nonlinear time-dependent PDE-constrained optimization problems. It is composed of a sequential quadratic programming (SQP) scheme to solve the PDE-constrained problem in an offline phase, a proper orthogonal decomposition (POD) approach to identify a lower dimensional solution space, and a neural network (NN) for fast online evaluations. The proposed method is showcased on a regularized least-square optimal control problem for the viscous Burgers' equation. It is concluded that significant online speed-up is achieved, compared to conventional methods using SQP and finite elements, at a cost of a prolonged offline phase and reduced accuracy. 

%% file: intro.tex
\section{Introduction}\label{sec:intro}
Nonlinear model predictive control is a well-established technique to solve optimal control problems within science, engineering and other areas \cite{qin2003survey, camacho2007nonlinear}. However, in cases where the underlying mathematical model is of high dimension, as is the case for PDEs, the optimization requires substantial computational power and time and cannot meet real-time constraints. In such mathematical models the curse of dimensionality is often a bottleneck \cite{bellman2015adaptive}. As a consequence, black box optimizers are not sufficient and one must resort to customized solvers. The state of the art customized solvers often rely on preconditioned iterative methods \cite{christiansen2018new, pearson2013fast}, high-performance computing \cite{biros2005parallel} and/or reduced order modeling techniques \cite{haasdonk2017reduced,biegler2003large} to improve the online computation time.\newline \indent
The contribution of this paper towards real-time optimization of large-scale processes, is a scheme based on proper orthogonal decomposition (POD), sequential quadratic programming (SQP) and neural networks (NNs), denoted the POD-SQP-NN scheme. This approach is in the field of reduced order modeling. The work is inspired by \cite{hesthaven2018non}, where a similar approach is presented for forward solution of nonlinear PDEs. The idea is to identify a low dimensional representation of the optimal control function, and then train a neural network to map inputs to this representation. The scheme utilizes a data driven approach to compute the low dimensional representation. One starts by exploring the solution space by computing optimal control functions, using the SQP algorithm \cite{troltzsch2010optimal}, and then computes a low dimensional representation using POD \cite{quarteroni2015reduced}. The idea is closely related to that of principal component analysis for dimensionality reduction of data representations \cite{friedman2001elements}. This approach exploits the often existing low dimensional structures to decrease the required computational effort and increase the accuracy of an NN. \newline \indent
The computation of the low dimensional representation and the training of the NN is kept in the offline stage. Thus, it is only necessary to evaluate the NN for given input in the online stage, which is an inexpensive procedure. This way one circumvents the problem associated with the curse of dimensionality in the online stage by keeping the full problem in the offline stage. \newline \indent
The POD-SQP-NN scheme is tested on a regularized least-square optimal control problem for the viscous Burgers' equation to demonstrate its potential. The viscous Burgers' equation has a nonlinear advection term. When this term is dominating shock waves will appear in the uncontrolled case, which makes it challenging to control \cite{troltzsch2001sqp}. \newline \indent
The paper is organized as follows. Section 2 introduces the optimal control problem and the optimality system. Section 3 presents the reduced order modeling framework. Section 4 outlines the use of NNs. Section 5 presents numerical results and conclusions, as well as future prospects, are made in section 6 and 7.

%% file: theory.tex
\section{Optimal Control of Nonlinear PDEs}\label{sec:optimal_control} 
Consider the class of nonlinear time-dependent PDEs of the form
\begin{subequations}
\begin{alignat}{2}
&\partial_t y - \epsilon \partial_{xx}y + N(y) = f,&\quad  &\text{in} \quad Q \\
&y(x,t) = 0, & &\text{on} \quad \Sigma   \\
&y(x,0) = y_0, & &\text{in} \quad \Omega  
\end{alignat}\label{PDE}
\end{subequations}
where $Q=\Omega\times (0,T)$ is a space-time cylinder, $\Omega\in\mathbb{R}^d$, $d=1,2,3$, $\Sigma=\Gamma\times (0,T)$, $\Gamma=\partial \Omega$, $\epsilon>0$, $f\in L^\infty(Q)$ is an external forcing, and $N$ is nonlinear operator. The initial value problem (IVP), \eqref{PDE}, is a general form of advection-reaction-diffusion equations with linear diffusion. Processes of this type encapsulates a great amount of physical phenomena, such as chemical reactions \cite{christiansen2018new}, fluid flows \cite{kuzmin2015finite}, and predator-prey systems \cite{ainseba2008reaction}. 

\subsection{The Optimal Control Problem}
To control \eqref{PDE} towards a desired state we formulate the PDE-constrained optimization problem:
\begin{subequations}
\begin{alignat}{2}
\min_{y,u} \quad &J(y,u), & &  \\
\text{s.t.} \quad &\partial_t y - \epsilon \partial_{xx}y + N(y) = f + bu,&\quad  &\text{in} \quad Q \\
&y(x,t) = 0, & &\text{on} \quad \Sigma   \\
&y(x,0) = y_0, & &\text{in} \quad \Omega  
\end{alignat}\label{optimal_control}
\end{subequations}
The goal of problem \eqref{optimal_control} is to determine the optimal control function, $u\in L^2(Q)$, such that the state, $y$, is as close as possible to a pre-defined desired state, $y_d\in L^2(Q)$. For this purpose a typical choice of the cost functional, $J$, is given by the tracking-type functional:
\begin{align}
J(y,u) = \frac{1}{2} \left|\left| y-y_d\right|\right|_2^2 + \frac{\beta}{2}  \left|\left| u\right|\right|_2^2, \label{cost_functional}
\end{align}
where $||\cdot||_2$ denotes the standard $L^2(Q)$ norm and $\beta>0$ is the Tikhonov regularization parameter, ensuring that the optimal control has higher regularity \cite{troltzsch2010optimal}. The control can be restricted to a subdomain, $\Omega_c\subset \Omega$ by choosing $b=b(x,t)$ as a characteristic function, i.e.  $b=b(x,t)\chi_{\Omega_c}(x,t)$. Note that in many real-world applications \eqref{optimal_control} is often accompanied by constraints on the control function, $u$. However, to illustrate the ideas and principles of this paper we leave out such constraints on $u$. 

\subsection{Optimality System}
Solving \eqref{optimal_control} can be done in more than one way. This paper follows the \textit{optimize-then-discretize} approach \cite{yilmaz2014all}. That is, we start by deriving the optimality conditions on the continuous level and then we discretize the resulting system. Utilizing the formal Lagrange method \cite{troltzsch2010optimal}, we arrive at the following system:
\begin{subequations}
\begin{alignat}{2}
&\partial_ty - \epsilon \partial_{xx} y + N(y) -\frac{b^2}{\beta}p= f, &\quad &\text{in} \quad Q, \label{first_order_optimality_general_problem_state1_reduced}    \\
-&\partial_tp - \epsilon \partial_{xx} p + N_y(y)p = y-y_d, & &\text{in} \quad Q, \label{first_order_optimality_general_problem_adjoint1_reduced}    \\
&y(x,t) = p(x,t) = 0, & &\text{on} \quad \Sigma,  \label{first_order_optimality_general_problem_state2_reduced}   \\
&y(x,0) = y_0, & &\text{in} \quad \Omega, \label{first_order_optimality_general_problem_state3_reduced}  \\
&p(x,T) = 0, & &\text{in} \quad \Omega. \label{first_order_optimality_general_problem_adjoint3_reduced} 
\end{alignat} \label{optimality_system}
\end{subequations}
where $p$ is the adjoint state. The control is recovered by $u=\frac{b}{\beta}p$. The IVP \eqref{optimality_system} is a system of nonlinear PDEs with the state evolving forward in time and the adjoint state backwards. We use a Newton/SQP approach to overcome the nonlinearities, which first linearizes \eqref{optimality_system} and then solve the resulting system until a given stopping criteria is fulfilled \cite{troltzsch2010optimal,troltzsch2001sqp}. Denoting by $(y_i,p_i)$, $1\leq i\leq  k$, the Newton iterates and $(y_{k+1},p_{k+1}):=(y,p)$, the next iterate, we recover the following linear system to solve for $(y,p)$:
\begin{subequations}
\begin{align}
&\partial_t y -\epsilon \partial_{xx} y + N_y(y_k) y - \frac{b^2}{\beta} p =F_1,  \\
-&\partial_t p -\epsilon \partial_{xx} p + N_{y}(y_k)p-y+N_{yy}(y_k)p_k y = F_2,
\end{align}\label{linear_SQP_system}
\end{subequations}
where
\begin{subequations}
\begin{align}
&F_1 := f + N(y_k)-N_y(y_k)y_k, \\
&F_2 := -y_d + N_{yy}(y_k)p_ky_k,
\end{align}
\end{subequations}
and with initial and boundary conditions as given in \eqref{optimality_system}. \newline \indent
While \eqref{linear_SQP_system} are linear, the problem of opposite time evolutions of the state and the adjoint state remains. To solve this problem we solve for all time steps at once. This is also referred to as an \textit{all-at-once} approach \cite{yilmaz2014all}.\newline \indent
Using the Galerkin finite element method (FEM), we aim to represent the solutions by
\begin{align}
y = \sum_{i=1}^{N_\delta} \hat{y}_i(t) \phi_i(x), \quad  p = \sum_{i=1}^{N_\delta} \hat{p}_i(t) \phi_i(x),
\end{align}
where $\phi_i\in H^1(\Omega)$ are the finite element basis functions and $\hat{y}_i, \hat{p}_i$ are the generalized Fourier coefficients, i.e. the degrees of freedom, to be estimated. We denote the finite dimensional solution space $V_\delta = \text{span}\left\{\phi_i\right\}_{i=1}^{N_\delta}$, introduce a time discretization, $t_i = i\Delta t$, $i=1\ldots N_t$, and define
\begin{subequations}
\begin{align}
&\textbf{y} = (\hat{y}_1(t_1),\ldots,\hat{y}_1(t_{N_t}),\ldots,\hat{y}_{N_\delta}(t_1),\ldots,\hat{y}_{N_\delta}(t_{N_t})),\\
&\textbf{p} = (\hat{p}_1(t_1),\ldots,\hat{p}_1(t_{N_t}),\ldots,\hat{p}_{N_\delta}(t_1),\ldots,\hat{p}_{N_\delta}(t_{N_t})),
\end{align}
\end{subequations}
as the vectors collecting the degrees of freedom for all time steps. Thus, we solve a system of the form
\begin{align}
\underbrace{\begin{bmatrix}
A_k & B_k^T \\ B_k^T & -C_k
\end{bmatrix}}_{L_k}
\begin{bmatrix}
\textbf{y} \\ \textbf{p}
\end{bmatrix}
=
\begin{bmatrix}
\textbf{F}_1 \\ \textbf{F}_2
\end{bmatrix}, \label{linear_matrix_system}
\end{align}
in each Newton iteration. Note that $L_k\in \mathbb{R}^{2N_\delta N_t}$ can potentially be a large matrix. Thus, solving \eqref{linear_matrix_system} multiple times may be expensive, which motivates the use of specialized approaches to improve computation time.

\section{Reduced Order Modeling}\label{sec:reduced_models}
Decreasing the size of \eqref{linear_matrix_system} is one way to reduce the computational cost. To do so we need an efficient reduced basis, i.e. a basis comprising fewer functions without compromising the accuracy for given set of parameters. Hence, we seek solutions of the form
\begin{align}
u_\texttt{rb} = \sum_{i=1}^{N_\texttt{rb}} \hat{u}_i(t) \psi_i(x),
\end{align}
where $N_\texttt{rb}\ll N_\delta$ and $\psi_i$ is the reduced basis functions. We denote the reduced space $V_\texttt{rb} = \text{span}\left\{\psi_i\right\}_{i=1}^{N_\texttt{rb}}\subset V_\delta$. This also motivates naming $V_\delta$ the \textit{full-order} or \textit{high-fidelity} space. To derive such a basis, we utilize the POD \cite{quarteroni2015reduced,hesthaven2016certified}. 

\subsection{Parametrized PDEs}
To describe the POD method, we begin by introducing the concept of parametrized PDEs. Here, we assume that the solution does not only depend on space and time but also on the parameters. Define the vector, $\boldsymbol{\mu} \in \mathbb{P} \subset \mathbb{R}^P$, to be the vector of all parameters. Such parameters could be either physical, such as the diffusion rate, characteristics of the initial condition, or geometric, such the shape and size of the domain. The input parameters will often comprise an initial state in optimal control problems. Therefore, it makes sense to decompose the parameter vector into an initial state and other parameters, $\boldsymbol{\mu} = (\textbf{y}_0,\boldsymbol{\lambda})$. \newline \indent
The set of all solutions for varying parameters makes up a manifold, denoted the \textit{solution manifold}:
\begin{align}
\mathcal{M} = \left\{u(\boldsymbol{\mu}) \: | \: \boldsymbol{\mu} \in \mathbb{P} \right\}.
\end{align} 
Similarly, we define $\mathcal{M}_\delta$ and $\mathcal{M}_\texttt{rb}$ to be the solutions manifold for the high-fidelity and the reduced basis approximations, respectively.

\subsection{Proper Orthogonal Decomposition}
The POD method to compute the reduced basis is a data driven approach. One starts by computing a series of  \textit{snapshots}. A snapshot is a solution of the full problem for a given set of parameters and time instance. Ideally the snapshots should make up a good representation of the solution manifold. Then we seek to minimize
\begin{align}
\sqrt{\sum_{i=1}^{N_s}\inf_{v\in V_\texttt{rb}}\left|\left|u_\delta^i- v\right|\right|_{L^2(\Omega)}^2 },
\end{align}
where $V_\texttt{rb}$ is the reduced space, $u_\delta^i\in V_\delta$ denotes a snapshot, i.e. the solution for a specific parameter choice and time instance, and $N_s=N_pN_t$ is the number of snapshots, where $N_p$ and $N_t$ are the number of parameter samples and time steps respectively. According to the orthogonal projection theorem the minimizing reduced basis representation of $u_\delta^i$ is the orthogonal projection\cite{quarteroni2015reduced}:
\begin{align}
\sqrt{\sum_{i=1}^{N_s} \left|\left|u_\delta^i - \sum_{k=1}^{N_\texttt{rb}}\langle u_\delta^i,\psi_j \rangle \psi_j \right|\right|_{L^2(\Omega)^2}}. \label{POD_min}
\end{align}
In the $l^2(\Omega)$ case one can show that minimizing \eqref{POD_min} is equivalent to the constrained maximization problem \cite{gubisch2017proper}:
\begin{subequations}
	\begin{align}
	\max_{\{\pmb{\psi}_i\}_{i=1}^N} \quad &\sum_{i=1}^{N_s} \sum_{j=1}^{N_\delta} \left|\textbf{u}^T_i \pmb{\psi}_j \right|^2, \\
	\text{s.t.} \quad & \pmb{\psi}_j^T\pmb{\psi}_i = \delta_{ij}, \quad \forall i,j
	\end{align}\label{max_reduced_problem}
\end{subequations}
where $\pmb{\psi}_j$ is the vector of coefficients of the reduced basis functions, expanded in terms of the FE functions, and $\textbf{y}_i$ is the degrees of freedom of a snapshot. By collecting all the snapshots in a snapshot matrix, $S$, and all the reduced basis functions as columns in a matrix, $W$, the solution to \eqref{max_reduced_problem} is the solution to the eigenvalue problem \cite{quarteroni2015reduced}:
\begin{align}
SS^T W = W\Lambda,
\end{align}
where $\Lambda$ is a diagonal matrix with the eigenvalues, corresponding to the columns of $W$, on the diagonal. The matrix $SS^T$ is denoted the correlation matrix \cite{quarteroni2015reduced}. Truncating the number of eigenvectors gives us the number of reduced basis functions in use. One transforms between the degrees of freedom, in the high-fidelity basis and the POD basis, by \cite{quarteroni2015reduced}:
\begin{align}
\textbf{y}_\delta = W \textbf{u}_\texttt{rb}, \quad \textbf{u}_\texttt{rb} = W^T \textbf{u}_\delta. \label{transform}
\end{align}
The projection error on the snapshots corresponds to the eigenvalues of the correlation matrix of the POD modes left out \cite{quarteroni2015reduced}:
\begin{align}
\sum_{i=1}^{N_s}||\textbf{u}_i -\underbrace{ WW^T \textbf{u}_i}_{\text{Projection}} ||_2^2 = \sum_{i = N_\texttt{rb}+1}^r \lambda_i.
\end{align}
Hence, the eigenvectors corresponding to the largest eigenvalues contain the most information about the solution space. Thus, one should choose the reduced basis functions by looking at the size of their respective eigenvalues.

\section{Artificial Neural Networks}\label{sec:neural_networks}
First, we will give a brief recap of the fundamental concepts regarding NNs. An NN can be considered a map, $F:\mathbb{R}^{M_I}\rightarrow \mathbb{R}^{M_O}$, that maps a vector of $M_I$ features to a vector of $M_O$ responses. In the case of a feedforward network, also denoted a perceptron, the map consists of a series of function compositions \cite{Wolf18mathfound}:
\begin{align}
F(x;\theta,H,L) = u_{L+1} \circ f_{L} \circ u_{L} \circ \ldots \circ f_1\circ u_1 (x), 
\end{align}
where $L$ is the number of hidden layers, $H$ is the number of neurons in each layer, $\theta$ is the set of weights and biases, $f_i$ are the activation functions, and $u_i$ are the propagation functions. The evaluation of a NN is called a forward propagation since it can be considered a propagation of the input features through the series of functions. The propagation function is often a linear combination with a bias term:
\begin{align}
u_i(x) = w_{i-1} x + \beta_i,
\end{align}
where $w_i$'s are matrices containing the weights and $\beta_i$'s are vectors containing the bias terms.
There are many choices for the activation functions. Recently the functions called rectified linear units (ReLUs) have gained much success due to great results. In this paper we make use of the slightly modified version called the leaky ReLU \cite{maas2013rectifier}:
\begin{align}
f(x) = \max(\alpha x, x),
\end{align}
where $\alpha$ is a small positive real number.\newline \indent
The motivation for using neural networks are their versatility. It can be shown that any continuous function $\varphi\in \mathbb{R}^{d}\rightarrow \mathbb{R}^{d'}$ can be approximated arbitrarily well by a deep neural network \cite{Wolf18mathfound}. However, in practice it is difficult to get arbitrary precision due to various problems. These include i) a large parameter space, ii) nonlinearities, and iii) a non-convex optimization problem for determining the parameters. 
\newline \indent Much time has gone into developing optimization algorithms for computing the weights and biases to get as close as possible to the theoretical potential. Second order methods can not, in general, be used since they only ensure local convergence and are too expensive in high dimensions. Therefore, stochastic variations of gradient descent algorithms, such as Adam, are often used \cite{kingma2014adam}.

\section{The POD-SQP-NN Scheme}
\begin{figure*}
	\centering
	\includegraphics[width=0.85\textwidth]{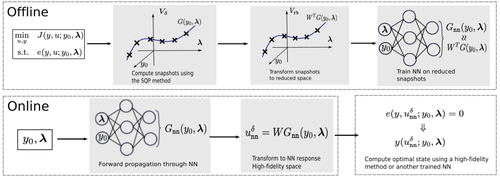}
	\caption{Workflow for the POD-SQP-NN scheme divided into the offline and online phase. Note the significantly heavier computations in the offline phase in shape of repeated solving by the SQP method and the training of a neural network.}
	\label{fig:procedure}
\end{figure*}

For PDE-constrained optimization problems one can consider a map, that takes input parameters and maps them to the high-fidelity optimal control solution manifold:
\begin{align}
G:\mathbb{P}\subset\mathbb{R}^P\rightarrow \mathbb{R}^{N_\delta N_t}, \quad (y_0,\boldsymbol{\lambda}) \mapsto G(y_0,\boldsymbol{\lambda}) = \textbf{u}_\delta.
\end{align}
$G(y_0,\boldsymbol{\lambda})$ is computed using the SQP method and finite elements. The image of the mapping, $G$, corresponds to the high-fidelity solution manifold for the control function. The general idea of the POD-SQP-NN scheme is to approximate this map in the POD space by using NNs. That is, by using the relation in \eqref{transform} we aim to compute a map, $G_\texttt{nn}$, such that
\begin{align}
G_\texttt{nn}(y_0,\boldsymbol{\lambda}) \approx W^T G(y_0,\boldsymbol{\lambda}), \quad \forall \boldsymbol{\mu} \in \mathbb{P},
\end{align}
from which we recover the optimal control in the high-fidelity space by 
\begin{align}
\textbf{u}_\texttt{nn}^\delta = W G_\texttt{nn}(y_0,\boldsymbol{\lambda}).
\end{align}
Consequently, the POD-SQP-NN scheme is decomposed into two stages:
\begin{itemize}
\item \textbf{The offline stage.}
\begin{enumerate}
\item Compute snapshots, $\left\{(y_0,\boldsymbol{\lambda})_i,G((y_0,\boldsymbol{\lambda})_i)\right\}_{i=1}^{N_s}$, using SQP and finite elements.
\item Compute POD basis functions and transform snapshots to get training data: $\left\{(y_0,\boldsymbol{\lambda})_i,W^TG((y_0,\boldsymbol{\lambda})_i)\right\}_{i=1}^{N_s}$ .
\item Train the neural network, $G_\texttt{nn}$, on training data. 
\end{enumerate}
\item \textbf{The online stage.}
\begin{enumerate}
\item Input parameters, $(y_0,\boldsymbol{\lambda})$, and compute $G_\texttt{nn}(y_0,\boldsymbol{\lambda})$ by forward propagation in NN.
\item Transform output to high-fidelity space, $\textbf{u}_\texttt{nn}^\delta = WG_\texttt{nn}(y_0,\boldsymbol{\lambda})$. 
\item (Optional) Compute the optimal state by a forward solving of the state equation.
\end{enumerate}
\end{itemize}
This procedure is outlined in Fig. \ref{fig:procedure}. It is apparent that the offline stage is quite extensive. It consists of several steps which are all, by themselves, computationally heavy tasks and are highly dependent on the dimension of the problem. However, the online stage is not computationally demanding, which is essential for real-time applications.  \newline \indent
By approximating the solution in the reduced space we reduce the number of output neurons and, thereby, also the number of weights to be trained significantly. This reduces training time as well as online computation time. Furthermore, it increases the accuracy of the NN since it reduces the size of the parameter space. \newline \indent
A typical problem when training neural networks is lack of enough data. In the proposed method, as much data as needed can be generated, given no time constraints on the offline phase. This, however, does not ensure high precision due to other constraints in neural networks such as  bias-variance trade-off, irreducible errors, etc. \cite{friedman2001elements}. However, it is worth discussing potential problems regarding the dimension of the input. Due to the curse of dimensionality it might be infeasible to compute the proposed map $G_\texttt{nn}$ since it might require an insurmountable amount of simulations to explore the full solution manifold. However, in many cases it can be exploited that the initial state lies on a low dimensional surface parametrized by only a few parameters, which reduces this problem significantly.\newline \indent
We define $\Pi_{\texttt{rb}}$ to be the projection operator onto the reduced basis space. By using the triangle inequality, one gets rough a priori error estimate:
\begin{align}\label{NN_approx_error}
\begin{split}
\left|\left|u-u_{\texttt{nn}}^\delta\right|\right|_{L^2(Q)} \leq &\underbrace{\left|\left|u- u_\delta\right|\right|_{L^2(Q)}}_{=: \varepsilon_\delta} + \underbrace{\left|\left|u_\delta - \Pi_{\texttt{rb}}u_\delta\right|\right|_{L^2(Q)}}_{=:\varepsilon_\texttt{rb}} \\ 
&+ \underbrace{\left|\left|\Pi_{\texttt{rb}}u_\delta- u_{\texttt{nn}}^\delta \right|\right|_{L^2(Q)}}_{=:\varepsilon_{\texttt{nn}}},
\end{split}
\end{align}
$\varepsilon_\delta$ is the error between the high-fidelity and the exact solution. $\varepsilon_\texttt{rb}$ is the discrepancy between the high-fidelity solution and its POD basis representation. It depends on the quality of the POD method which, in turn, depends on the how reducible the problem is and how well the solution manifold is sampled. In many cases this error decreases exponentially \cite{hesthaven2016certified}. Lastly, $\varepsilon_{\texttt{nn}}$ is the discrepancy between the POD representation and the NN output. It depends on the quality of the training of the network and is often the most difficult to drive down through training.

%% file: ex.tex
\section{Numerical example}\label{sec:ex}
As a test example, we consider the viscous Burgers' equation as the state equation. Thus, the PDE-constrained optimization problem is:
\begin{subequations}
	\begin{align} 
	\min_{u,y},\quad &J(y,u), \\
	s.t. \quad &\partial_ty - \epsilon \partial_{xx}y+ y\partial_xy = f+u, \quad \text{in} \quad Q,  \label{burgers_state} \\
	&y(x,t) = 0, \quad \text{on} \quad \Sigma ,  \\
	&y(x,0) = y_0, \quad \text{in} \quad \Omega,
	\end{align} \label{burgers_control_problem}
\end{subequations} 
with $J$ defined as in \eqref{cost_functional} and $f\in L^\infty(Q)$ a known forcing on the system. The well-posedness of the \eqref{burgers_control_problem} is discussed in \cite{troltzsch2001sqp}. The resulting linear problem to be solved in each Newton iteration is
\begin{subequations}
	\begin{align}
	&\partial_t y -\epsilon \partial_{xx}y + \partial_x (y_k y) - \frac{1}{\beta} p =f + y_k\partial_x y_k, \\
	-&\partial_t p -\epsilon \partial_{xx}p - y_k\partial_xp +(\partial_xp_k-1)y = -y_d + y_k\partial_xp_k,
	\end{align}
\end{subequations}
For the spatial discretizations FEniCS is used \cite{AlnaesBlechta2015a} and for the NN the PyTorch framework is used \cite{paszke2017automatic}. The unconditionally stable implicit Euler is used for the time discretization. 

\subsection{Case Study}
\begin{figure}
	\centering
	\includegraphics[width=0.35\textwidth]{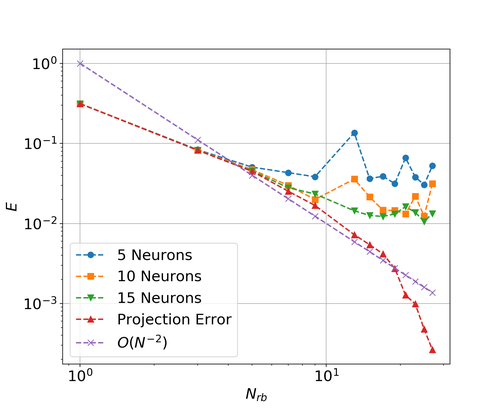}
	\caption{Convergence of the POD-SQP-NN method for varying number of neurons in each layer.}
	\label{fig:NN_convergence_adjoint}
\end{figure}

\begin{figure}
	\centering
	\includegraphics[width=0.35\textwidth]{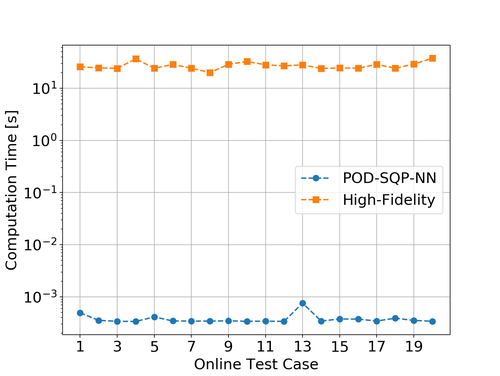}
	\caption{Execution time in seconds.}
	\label{fig:NN_time_burgers_case1}
\end{figure}

\begin{figure}
	\centering
	\includegraphics[width=0.35\textwidth]{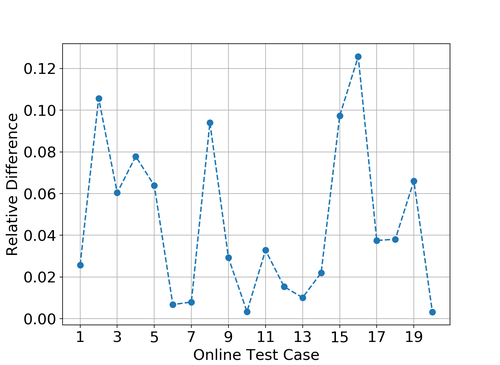}
	\caption{Relative difference in objective function between the high-fidelity and the POD-SQP-NN solution.}
	\label{fig:rel_diff}
\end{figure}

For the specific case study we choose $\Omega=[0,1]$, $T=0.5$ and the initial state:
\begin{align}
y_0(x) = 
\begin{cases}
h, & \text{for  } \omega \leq x \leq 1-\omega \\
0 & \text{otherwise}
\end{cases}.
\end{align}
The desired state is to retain the initial condition over time. This problem is quite similar to the one considered in \cite{troltzsch2001sqp}. We define the vector of parameters to be varied as $\boldsymbol{\mu}=(h,\omega,\epsilon)\in  [0.5,1] \times [0.1,0.3]\times [0.1,0.001]$. It is worth noting that $h$ and $\omega$ parametrizes the initial condition while $\epsilon$ defines the diffusion rate. To sample the parameters the Latin Hypercube Sampling method is utilized to ensure a random distribution and that the full parameter space is explored \cite{cochran2007sampling}.  \newline \indent
In this paper, we use $N_p=300$ and $N_t=75$ which in total gives $22500$ snapshots to compute the POD basis. Note that only 300 SQP all-at-once method solutions are needed to get all the snapshots. It is likely that that fewer snapshots would be sufficient since $f$ is static and $y_d=y_0$. \newline \indent
To train the neural network in the POD-SQP-NN scheme a training set of 2500 optimal control solutions is used. For the architecture we train a feedforward network with three hidden layers with $H$ neurons in each. Furthermore, the Adam optimizer was utilized using mini-batches of size 32 together with early stopping and the mean squared error with $L^2$-regularization as the loss function. Furthermore, the weights were initialized according to the standard Gaussian distribution. \newline \indent
The error is measured to be the relative difference of the neural network response and the high-fidelity solution:
\begin{align}
E = \frac{\left|\left| \textbf{u}_\texttt{nn}^\delta - \textbf{u}_\delta  \right|\right|_{L^2}}{\left|\left|\textbf{u}_\delta\right|\right|}.
\end{align}
From Fig. \ref{fig:NN_convergence_adjoint} it is clear that the accuracy is increasing with the number of POD basis functions. Furthermore, more neurons in each layer gives a better approximation until a certain point at which the error stagnates. However, it is worth noting that the error is not decreasing monotonically, but is rather noisy. This points to the fact that the networks are not trained well enough in every configuration. In the best cases we reach an accuracy in the order of magnitude $10^{-2}$ seconds.\\
In Fig. \ref{fig:NN_time_burgers_case1} we see the online computation time for the high-fidelity as well as for the POD-SQP-NN scheme. The high-fidelity scheme takes between 10 and 100 seconds while the POD-SQP-NN scheme is of the order of $10^{-3}$ seconds. This does not come as a surprise, since the online phase is just a single forward propagation for the POD-SQP-NN scheme.\newline \indent
Whether an accuracy of $10^{-2}$ is sufficient is highly dependent on the problem at hand. However, for many problems it is enough. In Figure \ref{fig:rel_diff} the relative difference in the objective functional,
\begin{align}
\frac{\left|\left| J(\textbf{y}_\texttt{nn}^\delta,\textbf{u}_\texttt{nn}^\delta) - J(\textbf{y}_\delta,\textbf{u}_\delta) \right|\right|_{L^2}}{\left|\left|J(\textbf{y}_\delta,\textbf{u}_\delta)\right|\right|},
\end{align}
is plotted. In all online test cases the discrepancy is no larger than 0.12, which suggests that an increase in accuracy will not yield significant changes in the objective functional. Whether such differences are acceptable depends on the problem at hand.  

%% file: concl.tex
\section{Conclusions}\label{sec:concl} 
A new and efficient scheme for online computations of nonlinear PDE-constrained optimization problems is proposed in the paper. We first showed that the problem at hand can be projected onto a lower dimensional manifold using the POD method, which lays the foundation for the POD-SQP-NN method. The scheme significantly decreases the online computation times for the case involving the viscous Burgers' equation. It did, however, come with a decrease in accuracy and a prolonged offline phase. Regarding the decrease in accuracy it is possible that it does not pose any practical limitations due to the minor discrepancies in the objective functions. For these reasons, we see a potential in this method to be used to solve problems that are currently intractable in real-time. 

\section{Future Prospects} 
What was not discussed in sufficient detail is how the method performs on very high dimensional problems. In many NMPC applications it is not necessarily possible to parametrize the initial conditions by a small set of parameters, as was done here. In such cases the entire initial state must be the input, which could potentially cause problems. Furthermore, the reliability is not addressed as no rigorous a posteriori error bounds exists. This was, however, not a problem in the examples given. \newline \indent
In ongoing work, we will, furthermore, explore different NN architectures such as residual networks, and reinforcement learning which has shown to be well-suited for time dependent data. Alternatively, instead of using NN other new approaches can be utilized to address the curse of dimensionality, such as spectral tensor train decompositions \cite{bigoni2016spectral}. 